\input amstex 
\documentstyle{amsppt}
\input bull-ppt
\keyedby{bull324/car}

\topmatter
\cvol{27}
\cvolyear{1992}
\cmonth{October}
\cyear{1992}
\cvolno{2}
\cpgs{284-287}
\title Is the boundary of a Siegel disk a Jordan curve? 
\endtitle
\author James T. Rogers, Jr.\endauthor
\shortauthor{J. T. Rogers, Jr.}
\address      Department of Mathematics, Tulane 
University, New Orleans,
Louisianna 70118\endaddress
\date November 21, 1991 and, in revised form, March 24, 
1992\enddate
\thanks This research was partially supported by a COR 
grant from Tulane 
University\endthanks
\abstract Bounded irreducible local Siegel disks include 
classical Siegel disks
of polynomials, bounded irreducible Siegel disks of 
rational and entire 
functions, and the examples of Herman and Moeckel. We show 
that there are only
two possibilities for the structure of the boundary of 
such a disk: either the 
boundary admits a nice decomposition onto a circle, or it 
is an indecomposable
continuum.\endabstract
\endtopmatter

\document

\heading 1. Introduction\endheading

Let $f\:\overline {\bold C}\rightarrow \overline {\bold 
C}$ be a rational map
of the Riemann sphere of degree at least two. The dynamics 
of $f$ divides $
\overline {\bold C}$ into two disjoint sets: the {\it 
stable set\/} or 
{\it Fatou set\/}, and the {\it unstable set\/} or {\it 
Julia set\/}. On the
Fatou set the dynamics of $f$ is well behaved, while the 
dynamics of $f$ on the
Julia set is chaotic.

The work of  Sullivan \cite{Su} completed the 
understanding of the dynamics of
$f$ on the  Fatou set. Every component of the Fatou set is 
eventually periodic,
and essentially five kinds of dynamical behavior are 
possible on these domains.
One of these behaviors is a Siegel disk.

A component $G$ of the Fatou set of $f$ is a {\it Siegel 
disk\/} if $f(G)=G$,
$G$ contains a neutral fixed point $w_0$, and $f\big|G$ is 
analytically
conjugate to a rotation. Siegel \cite S showed in 1942 
that such disks exist.

To say that $w_0$ is a {\it neutral fixed point\/} means 
$f(w_0)=w_0$ and
$|f\,'(w_0)|=1$. Hence $f\,'(w_0)=e^{2\pi i\theta}$ for 
some real number $
\theta$ in $[0,1)$. It is known that $\theta$ must be 
irrational, and much has
been written in the effort to decide which irrationals 
yield a Siegel disk (see 
\cite B or \cite M).

The dynamics of the Julia set is more subtle, and much is 
still unknown. Douady
and Sullivan \cite{D1} have raised a very natural 
question: Is the boundary of
a Siegel disk a Jordan curve? Herman \cite{D2} has 
obtained an affirmative
answer in special circumstances, but, in general, no 
answer is known.

More generally, let us define a {\it bounded local Siegel 
disk\/} to be  a pair
$(G,F_\theta)$, where $G$ is a bounded simply connected 
domain in $\bold C$,
and $F_\theta\:G\rightarrow G$ is a conformal map complex  
analytically
conjugate to a rotation through the irrational angle 
$\theta$ such that $F_
\theta$ extends continuously to the boundary of $G$. The 
fixed point $w_0$ is
again called a {\it Siegel point\/}. A bounded Siegel disk 
$(G,F_\theta)$
is {\it irreducible\/} if the boundary of $G$ separates 
the Siegel point $w_0$
from $\infty $, but no proper closed subset of the 
boundary has this property.
Bounded irreducible local Siegel disks include classical 
Siegel disks of 
polynomials as well as bounded irreducible Siegel disks of 
rational and entire 
functions and even the exotic examples of Herman \cite{H1} 
and Moeckel
\cite{Mo}.

We describe the structure of the boundaries of such 
domains by proving the 
following theorem.
\proclaim{Theorem 1.1} The boundary $\partial G$ of a 
bounded irreducible local 
Siegel disk satisfies exactly one of the following 
properties\RM:
\roster
\item The inverse $\varphi^{-1}$ of the Riemann map 
$\varphi\:\bold D 
\rightarrow G$ of the conjugation extends continuously to 
a map $\psi\:\partial
G\rightarrow\partial \bold D=S^1$, or
\item \<$\partial G$ is an indecomposable continuum.
\endroster
\endproclaim
An {\it indecomposable continuum\/} is a compact connected 
space which cannot
be written as a union $A\cup B$ with $A$ and $B$ connected 
closed proper
subsets of $X$. Indecomposable continua are complicated 
spaces; nevertheless,
Herman \cite{H1} has constructed a bounded irreducible  
local Siegel disk whose 
boundary is a certain  indecomposable continuum known as 
the pseudocircle.

In case (1), the point inverses of $\psi \:\partial 
G\rightarrow S^1$ are the
impressions of prime ends of $\varphi$. In particular, the 
point inverses of $
\psi$ are connected. A space with such a decomposition 
onto a circle can be
written as a union $A\cup B$ as described above, so the 
two possibilities are
mutually exclusive. Moeckel \cite{Mo} has constructed such 
an example in which
the point inverses of $\psi$ are either points or straight 
line intervals.

The Moeckel example shows that we cannot require 
$\varphi^{-1}$ to extend to a
homeomorphism in (1), while the Herman example shows that 
(2) can occur. Thus
the result is the best possible for such local Siegel disks.

The boundary of a Siegel disk is a Jordan curve if and 
only if the Riemann map
$\varphi\:\bold D \rightarrow G$ of the conjugation 
extends to a homeomorphism
of $\overline {\bold D}$ onto $\overline  G$. This is 
equivalent, of course, to
$\varphi^{-1}\:G\rightarrow \bold D$ extending to a 
homeomorphism of $
\overline  G$ onto $\overline {\bold D}$. Thus we may 
interpret the theorem to
imply any counterexample must be as nice as possible or as 
complicated as
possible.

The theorem above implies that a weak additional 
hypothesis is enough to answer
the Douady-Sullivan question affirmatively.

\proclaim{Theorem 1.2} Let $A$ be an arc in the boundary 
of a Siegel disk of a 
polynomial of degree $d\geq2$. If two internal rays from 
$G$ land on $A$, then
$\partial G$ is a Jordan curve.
\endproclaim

Thus, any arc in a counterexample must be ``hidden.'' In 
particular, we have
the following corollary.
\proclaim{Corollary 1.3} If the boundary of a Siegel disk 
of a polynomial of
degree $d\geq2$ is arcwise connected, then $\partial G$ is 
a Jordan
curve.\endproclaim

The Julia set of a polynomial $f$ is the closure of the 
set of repelling 
periodic points of $f$, and the boundary of a Siegel disk 
is a subset of the 
Julia set. Hence the next theorem is in one sense a little 
surprising.
\proclaim{Theorem 1.4} If the boundary $\partial G$ is a 
Siegel disk of a 
polynomial of degree $d\geq2$ contains a periodic point, 
then $\partial G$ is
an indecomposable continuum.
\endproclaim

This paper is an abstract of the results in \cite{R3}. The 
paper \cite{R3}
contains a brief history of indecomposable continua 
occurring in the study of
dynamical systems and suggests that it is not so 
unexpected that we must deal
with indecomposable continua in this situation.

\heading 2. The structure of the boundary of local Siegel 
disks\endheading

Let $(G,F_\theta)$ denote a bounded irreducible local 
Siegel disk. We need a
number of tools to complete the proof of the structure 
theorem.

The first is a result of Pommerenke and Rodin \cite{PR} 
about prime ends $\eta$
and their impressions $I(\eta)$.

\proclaim{Theorem 2.1} Each local Siegel disk has  a 
Pommerenke-Rodin
number\RM; i.e., there exists a number $d$ \RM(not to be 
confused with the
degree of a polynomial\RM) with $0\leq d\leq 2$ such that, 
for prime ends
$\eta_1$ and $\eta_2$ in $\partial \bold D$,
$$I(\eta_1)\cap I(\eta_2)\neq \emptyset \Leftrightarrow 
|\eta_1-\eta_2|\leq
d\.$$\endproclaim

The distance $|\eta_1-\eta_2|$ on $\partial \bold D$ is 
given by the Euclidean
metric on $\bold C$; hence, for example, 
$|\eta_1-\eta_2|=2$ if and only if $
\eta_1$ and $\eta_2$ are diametrically opposite. It 
follows that $d=0$ if and
only if all impressions are pairwise disjoint, while $d=2$ 
if and only if each
pair of impressions has a point in common.

The second and most important tool is the theory of prime 
ends as related to
indecomposable continua. The work of Rutt \cite{Ru} is 
used, for instance, in
proving the following result of the author \cite{R1}, a 
result that enables us
to recognize indecomposable continua by analytic methods.

\proclaim{Theorem 2.2} If $(G,\scr F_\theta)$ is a local 
Siegel disk, then $
\partial G$ is an indecomposable continuum if and only if 
there exists a prime
end $\eta$ of $G$ such that the impression 
$I(\eta)=\partial G$.\endproclaim

The proof of the structure theorem is completed by a 
somewhat delicate analysis
of the relationship between prime ends and indecomposable 
continua. The details
appear in \cite{R3}.

\enddocument